\documentclass{elsart}
\usepackage{mathrsfs}
\usepackage{amsfonts}

\usepackage{mathrsfs}

\usepackage{mathrsfs}
\usepackage{enumerate}
\usepackage{amssymb}
\usepackage{epsfig}
\usepackage{graphicx}
\usepackage{subfigure}
\usepackage{amsfonts}
\usepackage{amsmath}
\usepackage{amssymb}
\usepackage{latexsym,amscd}
\newcommand{\ucite}[1]{\textsuperscript{\textsuperscript{\cite{#1}}}}

\begin{document}

\begin{frontmatter}



\title{Toric arrangement and discrete truncated power}


\author{Renhong Wang},
\author{Mian Li\corauthref{cor1}}
\corauth[cor1]{Corresponding author.}\ead{limian3565@yahoo.cn}
\address{School of Mathematical Sciences, Dalian University of Technology, Dalian 116024, China}

\begin{abstract}
\hspace{0.4cm}  Discrete truncated power is very useful to study the
number of nonnegative integer solutions of linear Diophantine
equations.  In this paper, by using the Laplace transform and the
theory of toric arrangement, we show that discrete truncated power
is a periodic piecewise polynomial on the shifted integral lattice
cone. Based on the toric reduction method in the real field, we give
a toric arrangement method to compute discrete truncated power.
\end{abstract}


\end{frontmatter}
\section{Introduction}
\quad Let $\mathbf{Z}^{s}$ be the collection of $s$-dimensional
vectors whose components are integers, and $\mathbf{R}_{+}$ be the
collection of all nonnegative real numbers. Let $X$ be a multiset
(i.e. the set whose elements can be same) of $n$ vectors
$a_{1},\ldots,a_{n}$, where $a_i\in\mathbf{Z}^{s}\backslash\{0\}$,
$i=1,\ldots,n$, such that for $x_{1},\ldots,x_{n}\in\mathbf{R}_{+}$,
the equation $\sum_{i=1}^{n}x_{i}a_{i}=0$ does not have nonzero
solutions with respect to $a_i$, $i=1,\ldots,n$.  Thus we see that
\begin{equation}
t_{X}(x):=\#\{\beta=(\beta_{1},\ldots,\beta_{n})\in\mathbf{N}^{n}:
\sum_{i=1}^{n}\beta_{i}a_{i}=x\},
\end{equation}    is finite for every $x\in \mathbf{Z}^s$, where
$\#A$ denotes the cardinality of the set $A$,  and $\mathbf{N}^n$
denotes  the collection of $n$-dimensional vectors whose components
are  nonnegative integers. Then $t_{X}$ is called the discrete
truncated power or the partition function. Discrete truncated power
was first introduced by Dahmen and Micchelli \cite{Dahmen}.  They
exploit the relationship between discrete truncated power and
multivariate splines, and propose some important properties of
discrete truncated power \cite{Dahmen1}. In particular, they present
a recursive method to compute discrete truncated power, which is
based on the following recursive property:
$$t_{X}(\alpha)=\sum_{j=0}^{\infty}t_{X\backslash\{a_{i}\}}(\alpha-ja_{i}),~i=1,\ldots,n.$$
On the other hand, Laplace transform is found to be a useful tool to
study the discrete truncated power. To see this, we let
$$\Pi_{X}(\alpha):=\{x\in\mathbf{R}^{n}|\sum^{n}_{i=1}x_{i}a_{i}=\alpha,~x_{i}\geq 0\}$$
be a polytope. Then the value of $t_{X}(\alpha)$ is the number of
integral points in $\Pi_{X}(\alpha)$, and the Laplace transform of
discrete truncated power is given by:
\begin{equation}\label{tx}\sum_{h\in \mathbf{Z}^{s}}e^{-\langle h,x\rangle}t_{X}(h)
=\prod_{a\in X}\frac{1}{1-e^{-\langle a,x\rangle}}.\end{equation}
That motivates some scholars to study discrete truncated power in
another way. Szenes and Vergne establish the relation between
discrete truncated power and Jeffrey-Kirwan residues, and propose a
residue method to compute discrete truncated power \cite{Szenes}.
 Concini and Procesi present the relation between
discrete truncated power and toric arrangement, and show that
discrete truncated power is a quasi polynomial on some rational
lattice  \cite{Concini,Concini1}.

\quad In this paper, we investigate the discrete truncated power by
using toric arrangement. Our idea comes from two aspects: the
relation between multivariate spline and hyperplane arrangement
\cite{Concini}, and the Laplace transform of the discrete truncated
power. Although the second technique was used in
\cite{Szenes,Concini}, within the complex field, our argument and
result are simple and clear, as the discussion is based on the real
field. Moreover, we prove that the discrete truncated power $t_{X}$
is a periodic piecewise polynomial on the integral lattice cone
spanned by the vectors in $X$. This is more precise than the result
of Concini and Procesi in \cite{Concini}\footnote{In \cite{Concini},
the authors show that discrete truncated power is a quasi polynomial
on the rational lattice. However, as the discrete truncated power is
defined on the integral lattice, their result essentially makes
sense merely for the integral lattice. Hence our result can be
considered as an improvement of theirs.}. Moreover, we propose a
toric arrangement method to compute the discrete truncated power,
which is able to give the explicit expression of $t_X$.

\quad  The paper is organized as follows. In Section 2, we briefly
review the conclusion given by Concini and Procesi, and introduce
some useful definitions and notations. In Section 3, we propose a
new method to show the relation between discrete truncated power and
toric arrangement, by using factor decompositions in the real field.
We can see that our method can be used to compute  the explicit
expression of the discrete truncated power, which is also
illustrated in two examples.

\section{Toric arrangement and discrete truncated power}

\quad In this section, we introduce some definitions and notations,
and previous results regarding the relation between toric
arrangement and discrete truncated power. We make the stipulation
that in this paper, the set associated with the discrete truncated
power  is multiset.

\quad Let $\Lambda$ be an integral lattice in $\mathbf{R}^{s}$ and
$$C[\Lambda]:=\{e^{\langle a,x\rangle},~a\in \Lambda\}$$
 be a collection of multivariate exponential functions. Let $v$ be
an integral vector in $\mathbf{R}^{s}$ and $\tau_{v}$ be a
translation operator such that
$$\tau_{v}f(x)=f(x+v)$$    for any function $f$ on
$\mathbf{R}^{s}$. Then we define
$$\tau[\Lambda]=\{\tau_{v},~v\in \Lambda\}$$
to be the collection of translation operators associated with
$\Lambda$.  Denote $\mathcal {L}$, $\mathcal {L}^{-1}$ to be the
operators of the Laplace transform and the inverse Laplace transform
respectively. Then the elements in $C[\Lambda]$ and $\tau[\Lambda]$
have the following relations:  \begin{align*}&\mathcal
{L}(\tau_{v})=e^{\langle v,x\rangle},\\
&\mathcal {L}^{-1}(e^{\langle v,x\rangle})=\tau_{v}.\end{align*}
 Let $x^{\alpha}=x_{1}^{\alpha_{1}}x_{2}^{\alpha_{2}}\cdots
x_{s}^{\alpha_{s}}$ be a monomial in $s$ variables, where
$\alpha_{1},\alpha_{2},\cdots,\alpha_{s}$ are nonnegative integers,
and denote  $$\mathbf{P}^{(s)}(x):=\{\sum_{\alpha}
a_{\alpha}x^{\alpha},~\alpha\in \mathbf{N}^{s},~a_{\alpha}\in
\mathbf{R}\}$$   to be the space of all the polynomials in $s$
variables. Let
$D^{\alpha}=D_{x_{1}}^{\alpha_{1}}D_{x_{2}}^{\alpha_{2}}\cdots
D_{x_{s}}^{\alpha_{s}}$ be a difference operators in $s$ variables,
and denote  $$\mathbf{P}^{(s)}(D):=\{\sum_{\alpha}
a_{\alpha}D^{\alpha},~\alpha\in \mathbf{N}^{s},~a_{\alpha}\in
\mathbf{R}\}$$  to be the space of all difference operators in $s$
variables. Similarly,  the elements in $\mathbf{P}^{(s)}(x)$ and
$\mathbf{P}^{(s)}(D)$ have the following relations:
\begin{align*}& \mathcal {L}(x^{\alpha})=(-D)^{\alpha},\\
& \mathcal {L}^{-1}(D^{\alpha})=(-x)^{\alpha}.
\end{align*}
Because the Laplace transform of $t_{X}$: Eq.(\ref{tx}),  has the
form $$\prod_{a\in X}\frac{1}{1-e^{-\langle a,x\rangle}},$$  the
following theorem is given.
\begin{thm}\label{laplace} \hspace{-2mm}{\rm \ucite{Concini}}
Let $$\mathcal {S}_{X}=C[\Lambda]\mathbf{P}^{(s)}(D)\prod_{a\in
X}(1-e^{-\langle a,x\rangle})^{-1},~~~\mathcal
{T}_{X}=\tau[\Lambda]\mathbf{P}^{(s)}(x)t_{X}$$   be two collections
of functions. Then under the Laplace transform, $\mathcal {T}_{X}$
is mapped isomorphically onto $\mathcal {S}_{X}$. \end{thm}
This
theorem means that for any element $\tau_{c}p(x)t_{X}(x)$ in
$\mathcal {T}_{X}$, there exists an element $e^{\langle
c,x\rangle}p(D)\prod_{a\in X}(1-e^{-ax})^{-1}$ in $\mathcal {S}_{X}$
such that
$$\mathcal{L}(\tau_{c}(p(x)t_{X}(x)))=p(x+c)t_{X}(x+c)=e^{\langle
c,x\rangle}p(D)\prod_{a\in X}(1-e^{-\langle a,x\rangle})^{-1}.$$
Contrarily, for any element $e^{\langle c,x\rangle}p(D)\prod_{a\in
X}(1-e^{-ax})^{-1}$ in $\mathcal {S}_{X}$, there exits an element
$\tau_{c}p(x)t_{X}(x)$ in $\mathcal {T}_{X}$ such that $$\mathcal
{L}^{-1}(e^{\langle c,x\rangle}p(D)\prod_{a\in X}(1-e^{-\langle
a,x\rangle})^{-1})=p(x+c)t_{X}(x+c)=\tau_{c}(p(x)t_{X}(x)).$$

To study the properties of $t_{X}$ on the lattice, Concini and
Procesi propose a toric reduction method to study the discrete
truncated power. The result is given in the following theorem.

\begin{thm}\label{concini}\hspace{-2mm}{\rm \ucite{Concini}}
Let $X=\{a_{1},a_{2},\ldots,a_{N}\}$ be a set of $N$ points on
$\Lambda$, and $\widetilde{X}$ be all the vectors in $X$ with
positive rational multiples. Then the function $\prod_{a\in
X}\frac{1}{1-e^{-\langle a,x\rangle}}$ can be written as a linear
combination of the form $\prod_{a\in A}\frac{1}{(1-e^{-\langle
a,x\rangle})^{h_a}}$ with constant coefficients, where $A$ is a
linearly independent set of elements in $\widetilde{X}$.
\end{thm}

By using the preceding theorem, Concini and Procesi give the
following theorem, which describes the general structure of $t_{X}$.

\begin{thm}\label{concini1}\hspace{-2mm}{\rm \ucite{Concini}}
Let $A$ be a linearly independent set of elements in
$\widetilde{X}$, such that the function $\prod_{a\in
X}\frac{1}{1-e^{-\langle a,x\rangle}}$ can be written as a linear
combination of the form $\prod_{a\in A}\frac{1}{(1-e^{-\langle
a,x\rangle})^{h_a}}$ with constant coefficients. On the closure of
the cone spanned by $A$, the partition function $t_{X}$ coincides
with a quasi polynomial for some rational lattice $\Lambda/n$.
\end{thm}

The idea of Concini and Procesi gives us a new understanding of
discrete truncated power. Based on this idea, we get our main result
in the next section (Theorem \ref{result}).

\section{Main results}

\hspace{0.4cm}We present our results in this section. First of all,
we give a useful lemma.


%
%

\quad Let $y_{i}$, $\alpha_{i}$, $i=1,\ldots,s$ be $2s$ functions
from $\mathbf{R}^{s}$ to $\mathbf{R}$, such that $y_{i}$ is a
nonzero function for each $i$. Then we have the following lemma.
\begin{lem}\label{lem2}
Let $y_{0}=\alpha_{1}y_{1}+\alpha_{2}y_{2}+\cdots+\alpha_{s}y_{s}$
be a function and  $h_{i}$ be a positive integer, $i=1,\ldots,s$.
Then $1/(y_{0}\prod\limits^{s}_{i=1}y_{i}^{h_{i}})$ can be written
as a linear combination of elements of the form
$f/\prod\limits^{s}_{i=0\atop i\neq k}y_{i}^{t_{i}}$, where $f$ is
the product of some $\alpha_{i}^{m_{i}}$, $1\leq k\leq s$, $t_{i}$,
$m_{i}$ are integers, and $\sum\limits^{s}_{i=0\atop i\neq
k}t_i=1+\sum\limits^{s}_{i=1}h_i$.
\end{lem}
\begin{pf}
Assume that $\alpha_{i}$ is a nonzero function. Then   we  get the
following algorithm to decompose
$1/(y_{0}\prod\limits^{s}_{i=1}y_{i}^{h_{i}})$:
\begin{equation}\frac{1}{y_{0}\prod\limits^{s}_{i=1}y_{i}^{h_{i}}}
=\frac{(\alpha_{1}y_{1}+\cdots+\alpha_{s}y_{s})}{y_{0}^{2}\prod\limits^{s}_{i=1}y_{i}^{h_{i}}}
=\sum\limits_{j}\frac{\alpha_{j}}{y_{0}^{2}\prod\limits^{s}_{i=1
\atop i\neq j}y_{i}^{h_{i}}\cdot y_{j}^{h_{j}-1}}.\end{equation}
Then  $1/(y_{0}\prod\limits^{s}_{i=1}y_{i}^{h_{i}})$ can be
decomposed into a sum of $s$ terms.   The numerator of each term is
$\alpha_{i}$. In the denominator, there exists $y_{i}$ whose power
decreases by one after decomposition, and the sum of powers of all
$y_i$ is also $1+\sum\limits^{s}_{i=1}h_i$. By using the algorithm
again, we have:
\begin{equation}\label{equation1}
\begin{split}\frac{\alpha_{j}}{y_{0}^{2}y_{j}^{h_{j}-1}\prod\limits^{s}_{i=1 \atop i\neq j}y_{i}^{h_{i}}}&=
\frac{\alpha_{j}(\alpha_{1}y_{1}+\cdots+\alpha_{s}y_{s})}{y_{0}^{3}y_{j}^{h_{j}-1}\prod\limits^{s}_{i=1
\atop i\neq j}y_{i}^{h_{i}}}\\
&=\frac{\alpha_{j}^{2}}{y_{0}^{3}y_{j}^{h_{j}-2}\prod\limits^{s}_{i=1
\atop i\neq j}y_{i}^{h_{i}}}+\sum\limits^{s}_{k=1 \atop k \neq
j}\frac{\alpha_{j}\alpha_{k}}{y_{0}^{3}y_{j}^{h_{j}-1}y_{k}^{h_{k}-1}\prod\limits^{s}_{i=1
\atop i\neq j,k}y_{i}^{h_{i}}}.\end{split}\end{equation}  If the
denominator of a term in the right hand side of Eq.(\ref{equation1})
 is the form of $\prod\limits^{s}_{i=0 \atop i\neq
k}y_{i}^{t_{i}}$, $k\in \{1,\ldots,s\}$,  we stop using the
algorithm on this term. Then we can express the numerator of this
factor in the form of $\prod\limits^{s}_{i=1}\alpha_{i}^{m_{i}}$ by
using Eq.(\ref{equation1}), where $m_i$ is a nonnegative integer.
Because $\sum_{i=1}^{s}h_{i}$ is a finite number,
$1/(y_{0}\prod\limits^{s}_{i=0}y_{i}^{h_{i}})$ can be decomposed
into a linear combination of elements of the form
$\prod\limits^{s}_{i=1}\alpha_{i}^{m_{i}}/\prod\limits^{s}_{i=0\atop
i\neq k}y_{i}^{t_{i}}$, where $\sum\limits^{s}_{i=0\atop i\neq
k}t_i=1+\sum\limits^{s}_{i=1}h_i$.

\quad If some $\alpha_{i}\equiv0$, then we denote
$\sigma:=\{i|\alpha_{i}=0\}$ to be the collection of indices such
that $\alpha_{i}\equiv0$. As
$$\frac{1}{y_{0}\prod\limits^{s}_{i=0}y_{i}^{h_{i}}}
=\frac{1}{\prod_{i\in
\sigma}y_{i}^{h_{i}}}\cdot\frac{1}{y_{0}\prod_{i\in
\{1,\ldots,s\}\backslash\sigma}y_{i}^{h_{i}}},$$ We consider the
decomposition of $1/(y_{0}\prod_{i\in
\{1,\ldots,s\}\backslash\sigma}y_{i}^{h_{i}})$. By using the similar
argument, we  can get the result. The proof is finished. $\square$
\end{pf}

\hspace{0.4cm}We can see the result of Lemma \ref{lem2} is similar
with the reduction of hyperplane arrangement. However, our result
generalizes $\alpha_{i}$ from a number to a function.

\hspace{0.4cm}   Let $X\subset\mathbf{Z}^{s}\setminus\{0\}$ be a set
of integer vectors in $\mathbf{R}^{s}$ such that $\dim X=s$. Denote
$\widetilde{X}$ to be the collection of nonnegative integral
multiples of the vectors in $X$. Then we can get the following
reduction theorem which is different from Theorem \ref{concini}. We
shall explain the difference later.

\begin{thm}\label{result1}The function $\prod_{a\in
X}\frac{1}{1-e^{-\langle a,x\rangle}}$ can be written as a
combination of elements of the form $\prod_{a\in A}\frac{e^{\langle
c,x\rangle}}{(1-e^{-\langle a,x\rangle})^{h_{a}}}$, with $h_{a}\in
\mathbf{N}$, such that $$\sum_{a\in X}h_a=\#X,$$ where $A$ is a
linearly independent set of $s$ elements in $\widetilde{X}$ and $c$
is a linear combination of elements in $X$.
\end{thm}

\begin{pf}
We use the mathematical induction on the number of elements in $X$.
If $X$ is linearly independent, the result is obvious. Otherwise, we
assume that the conclusion is valid for $X$ with $N-1$ elements,
where $N\geq s+1$.  Let
$$X=\{a_{1},a_{2},\ldots,a_{_N}\}$$ be a set of $N$ vectors.
As $X$ is not independent, there is an element, say $a_{_N}$, which
depends on the rest of the system
$$X'=\{a_{1},a_{2},\ldots,a_{_{N-1}}\}.$$
From the assumption, the function $\prod_{a\in
X'}\frac{1}{1-e^{-\langle a,x\rangle}}$ can be represented as a
linear combination of elements of the form
$$\prod_{a\in
A'}\frac{e^{\langle c',x\rangle}}{(1-e^{-\langle
a,x\rangle})^{h_{a}}},$$ where $h_{a}\in \mathbf{N}$, such that
$\sum_{a\in A'}h_a=\#X'$, $c'=\sum_{a\in X'} a$ (here $a$ can be
repeated) and $A'$ is a linearly independent set of $s$ elements in
$\widetilde{X'}$ where $\widetilde{X'}=\{na|a\in X',~n\in
\mathbf{N}\}$. Then for the function $\prod_{a\in
X}\frac{1}{1-e^{-\langle a,x\rangle}}$, we have   $$\prod_{a\in
X}\frac{1}{1-e^{-\langle a,x\rangle}}=\frac{1}{(1-e^{-\langle
a_{N},x\rangle})}\sum \prod_{a\in A'}\frac{e^{\langle
c',x\rangle}}{(1-e^{-\langle a,x\rangle})^{h_{a}}},$$  where $\sum$
denotes the summation of all the terms with the form $ \prod_{a\in
A'}\frac{e^{\langle c',x\rangle}}{(1-e^{-\langle
a,x\rangle})^{h_{a}}}$.
 Now let us
consider the factor
\begin{equation}\label{factor}\frac{1}{(1-e^{-\langle a_{_N},x\rangle})} \prod_{a\in
A'}\frac{e^{\langle c',x\rangle}}{(1-e^{-\langle
a,x\rangle})^{h_{a}}}.\end{equation} If the factor (\ref{factor})
can be written as the linear combination of the form $$\prod_{a\in
A}\frac{e^{\langle c,x\rangle}}{(1-e^{-\langle
a,x\rangle})^{h_{a}}}$$ where $c$ is a linear combination of
elements in $X$ and $A$ is a linearly independent set of $s$
elements in $\{na|a\in X\},~n\in \mathbf{N}\}$, then the conclusion
is valid. Without loss of generality, we let
$A'=\{a_{1},a_{2},\ldots,a_{s}\}$ (then $h_a=h_{a_{i}}$). And let
$h_i=h_{a_{i}}$ for simplicity. Because
$$A'\bigcup\{a_{_N}\}\subset X\subset \mathbf{Z}^{s},$$ there exists some
$m_{i}\in \mathbf{N}$, $i\in \{1,2,\ldots,s,N\}$, such that
\begin{equation}\label{equation2}m_{_N}a_{_N}=\sum\limits_{i\in\sigma}m_{i}a_{i}-\sum\limits_{i\in\{1,\ldots,s\}\backslash
\sigma}m_{i}a_{i},\end{equation}
where
$\sigma\subset\{1,\ldots,s\}$.

\quad Assume  $m_i\neq0$, $i=1,\ldots,s,N$. By reordering the
elements in $A$:
\begin{equation*}m_{_N}a_{_N}=\sum\limits_{i=1}^{k}m_{i}a_{i}-\sum\limits_{i=k+1}^{s}m_{i}a_{i},\end{equation*}
we have
\begin{eqnarray*}&&  1-e^{-\langle n_{_N}a_{_N},,x\rangle}\\
&=&1-e^{-\sum\limits_{i=1}^{k}\langle
n_{i}a_{i},x\rangle+\sum\limits_{i=k+1}^{s}\langle
n_{i}a_{i},x\rangle}=1-\prod_{i=1}^{k}e^{-\langle
n_{i}a_{i},x\rangle}\prod_{i=k+1}^{s}e^{\langle
n_{i}a_{i},x\rangle}\\
 &=&\prod_{i=1}^{k}e^{-\langle
n_{i}a_{i},x\rangle}\prod_{i=k+1}^{s-1}e^{\langle
n_{i}a_{i},x\rangle}(1-e^{\langle n_{s}a_{s},x\rangle})
+\cdots+\prod_{i=1}^{k}e^{-\langle n_{i}a_{i},x\rangle}(1-e^{\langle
n_{k+1}a_{k+1},x\rangle})\\
 &&+\prod_{i=1}^{k-1}e^{-\langle
n_{i}a_{i},x\rangle}(1-e^{-\langle
n_{k}a_{k},x\rangle})+\cdots+(1-e^{-\langle n_{1}a_{1},x\rangle})\\
&=&\prod_{i=1}^{k}e^{-\langle
n_{i}a_{i},x\rangle}\prod_{i=k+1}^{s-1}e^{\langle
n_{i}a_{i},x\rangle}(-e^{\langle
n_{s}a_{s},x\rangle})(\sum^{n_{s}-1}_{j=0}e^{-\langle
ja_{s},x\rangle})(1-e^{-\langle a_{s},x\rangle}) +\cdots\\
&&+\prod_{i=1}^{k}e^{-\langle n_{i}a_{i},x\rangle}(-e^{\langle
n_{k+1}a_{k+1},x\rangle})(\sum^{n_{k+1}-1}_{j=0}e^{-\langle
ja_{k+1},x\rangle})(1-e^{-\langle a_{k+1},x\rangle})\\
&&+\prod_{i=1}^{k-1}e^{-\langle
n_{i}a_{i},x\rangle}(\sum^{n_{k}-1}_{j=0}e^{-\langle
ja_{k},x\rangle})(1-e^{-\langle
a_{k},x\rangle})+\cdots+(\sum^{n_{1}-1}_{j=0}e^{-\langle
ja_{1},x\rangle})(1-e^{-\langle a_{1},x\rangle}).\end{eqnarray*}
To
give simple argument, we give the following notations:
$$\alpha_{t}=\left\{
   \begin{aligned}
   &1  &t=1,  \\
   &\prod_{i=1}^{t-1}e^{-\langle n_{i}a_{i},x\rangle}  &2\leq t\leq k+1,\\
   &\prod_{i=1}^{k}e^{-\langle
n_{i}a_{i},x\rangle}\prod_{i=k+1}^{t-1}e^{\langle
n_{i}a_{i},x\rangle}  &k+2\leq t\leq s,\\
   \end{aligned}
   \right.$$
$$\beta_{l}=\left\{
   \begin{aligned}
   &\sum_{i=0}^{l-1}e^{-\langle ia_{l},x\rangle}  &1\leq l\leq k,~l=N,  \\
   &-e^{\langle la_{l},x\rangle}\sum_{i=1}^{l-1}e^{-\langle ia_{l},x\rangle}  &k+1\leq l\leq s,\\
   \end{aligned}
   \right.~$$
$$y_{d}=\left\{
   \begin{aligned}
   &1-e^{-\langle a_{d},x\rangle}  &1\leq d\leq s,  \\
   &1-e^{-\langle n_{d}a_{d},x\rangle}  &d=N.\\
   \end{aligned}
   \right.~~~~~~~~~~~~~~~~~~~~$$
Then we have
$$y_{_N}=\sum_{i=1}^{s}\alpha_{i}y_{i}$$ and
$$\frac{1}{(1-e^{-\langle a_{_N},x\rangle})} \prod_{a\in
A'}\frac{e^{\langle c',x\rangle}}{(1-e^{-\langle
a,x\rangle})^{h_{a}}}=e^{\langle
c',x\rangle}\frac{\beta_{N}}{y_{_N}}\prod_{i=1}^{s}\frac{\alpha_{i}\beta_{i}}{y_{i}^{h_{i}}}=(e^{\langle
c',x\rangle}\beta_{N}\prod_{i=1}^{s}\alpha_{i}\beta_{i})\frac{1}{y_{_N}\prod_{i=1}^{s}y_{i}^{h_i}}.$$
According to Lemma \ref{lem2}, $1/(y_{_N}\prod_{i=1}^{s}y_{i})$ can
be represented as the linear combination of elements of the form
$\prod\limits^{s}_{i=1}\alpha_{i}^{n_{i}}/
(y_{_N}^{t_{N}}\prod\limits^{s}_{i=1\atop i\neq k}y_{i}^{t_{i}})$,
where $t_N+\sum\limits^{s}_{i=1\atop i\neq
k}t_i=1+\sum\limits^{i}_{i=1}h_i$. Because each $m_i$ is a nonzero
integer in Eq.(\ref{equation2}), $$\{y_i: ~i=1,\ldots
k-1,k+1,\ldots,s,N\}$$ are linearly independent. We observe that
\begin{equation}\label{e(cx)}(e^{\langle
c',x\rangle}\beta_{N}\prod_{i=1}^{s}\alpha_{i}\beta_{i})(\prod\limits^{s}_{i=1}\alpha_{i}^{n_{i}})
\end{equation}
is obviously a sum of the form of $e^{\langle c,x\rangle}$,
where $c$ is a linear combination of elements in $X$.

If some $m_i$'s are equal to 0 in Eq. (\ref{equation2}). We let
$$m_{_N}a_{_N}=\sum\limits_{i=1}^{d}m_{i}a_{i},$$ where integer $m_i\neq0$
for $i=1,\ldots d$. Then we consider the decomposition of
$$\frac{1}{(1-e^{-\langle a_{_N},x\rangle})} \prod_{i=1}^{d}\frac{e^{\langle c',x\rangle}}{(1-e^{-\langle
a_{i},x\rangle})^{h_{a_{i}}}}.$$ With the similar discussion, we
also can conclude that the hypotheses is right for
$X=\{a_1,\ldots,a_{_N}\}$. Then the proof is finished. $\square$
\end{pf}

\hspace{0.4cm}With Theorem \ref{laplace} and Theorem \ref{result1}
in hands, we can get the following theorem, which is  different from
but more precise than Theorem \ref{concini1}. Before proposing the
theorem, we introduce some definitions. For a matrix
$X=\{a_{1},\ldots,a_{N}\}\subset\mathbf{Z}^{s}$, we denote by
$$\Lambda_{X}^{+}:=\{\sum\alpha_{i}a_{i}|\alpha_{i}\in\mathbf{N}\}$$
a lattice cone and
$$\Lambda_{X}^{+}(c):=\{c+\sum\alpha_{i}a_{i}|\alpha_{i}\in\mathbf{N}\}$$
a shift of $\Lambda_{X}^{+}$. Let $C(X)$ be the cone
$$C(X):=\{\sum\alpha_{i}a_{i}|\alpha_{i}\in\mathbf{R}_{+}\}.$$
Then the following theorem is established.

\begin{thm}\label{result}
Let $X=\{a_{1},\ldots,a_{N}\}$ be a set of integral vectors in
$\mathbf{Z}^{s}$, then $t_{X}$ is a sum of some periodic piecewise
polynomials on the shifted sublattice cone of $\Lambda_{X}^{+}$.
\end{thm}
\begin{pf}
By Theorems \ref{laplace} and \ref{result1}, we have
$$\mathcal {L}(t_{X})=\frac{1}{\prod_{a\in X}(1-e^{-ax})}=\sum
\frac{e^{\langle c,x\rangle}}{\prod_{a\in A}(1-e^{-\langle
a,x\rangle})^{h_{a}}},$$ where $A=\{a_{1},\ldots,a_{s}\}$ is a basis
formed by vectors in $\widetilde{X}:=\{na|a\in X,~n\in\mathbf{N}\}$,
$c$ is a linear combination of elements in $X$ and $\sum_{a\in
A}h_a=\#X$. Let $h_i=h_{a_{i}}$ for simplicity. Now, let us compute
the inverse Laplace transform of $\frac{e^{\langle
c,x\rangle}}{\prod_{a\in A}(1-e^{-\langle a,x\rangle})^{h_{i}}}$.

Let $A_{i}^{\bot}$ be a vector who is perpendicular to all the
vectors in $A$ except for $a_{i}$. Then we have
\begin{equation}\label{equation}\frac{e^{\langle c,x\rangle}}{\prod_{a\in A}(1-e^{-\langle
a,x\rangle})^{h_{i}}}=e^{\langle
c,x\rangle}\prod_{i=1}^{s}\frac{(-1)^{h_{i}}(e^{\langle
a_{i},x\rangle}D_{A_{i}^{\bot}})^{h_{i}-1}}{(h_{i}-1)!\langle
A_{i}^{\bot},a_{i}\rangle^{h_{i}-1}}\cdot
\frac{1}{\prod_{i=1}^{s}(1-e^{-\langle a,x\rangle})},\end{equation}
where $D_{a}$ is a difference operator along vector $a$. Through
Eq.(\ref{equation}) and Theorem \ref{laplace}, we can show that the
inverse Laplace transform of $\frac{e^{\langle
c,x\rangle}}{\prod_{a\in A}(1-e^{-\langle a,x\rangle})^{h_{i}}}$ is
$$\tau_{c}\prod_{i=1}^{s}\frac{(\tau_{a_{i}}\langle
A_{i}^{\bot},x\rangle)^{h_{i}-1}}{(h_{i}-1)!\langle
A_{i}^{\bot},a_{i}\rangle^{h_{i}-1}}t_{A}=\prod_{i=1}^{s}\prod_{j=1}^{h_{i}-1}\frac{\langle
A_{i}^{\bot},x+ja_{i}+c\rangle}{(h_{i}-1)!\langle
A_{i}^{\bot},a_{i}\rangle^{h_{i}-1}}t_{A}(x+\sum_{k=1}^{s}
(h_{k}-1)a_{k}+c).$$ It is clearly that     $$t_{A}(x+\sum_{k=1}^{s}
(h_{k}-1)a_{k}+c)=t_{A}(x+c)\equiv1$$   on the lattice
$\Lambda_{A}^{+}(c)$. Thus the inverse Laplace transform of
$\frac{e^{\langle c,x\rangle}}{\prod_{a\in A}(1-e^{-\langle
a,x\rangle})^{h_{a}}}$ is a polynomial of degree $\#X-s$ on the
lattice $\Lambda_{A}^{+}(c)$. It is obvious that
$\Lambda_{A}^{+}\in C(X)$ is a sublattice cone of $\Lambda_{X}^{+}$.
The theorem is proved. $\square$
\end{pf}

\hspace{0.4cm}Now let us analysis the difference between Theorem
\ref{concini1} and Theorem \ref{result}. The main difference is
their domain, i.e. the lattice. For discrete truncated power
$t_{X}$, the lattice
$$\Lambda_{X}=\{\Sigma\alpha_{i}a_{i}|\alpha_{i}\in\mathbf{Z}\}$$ in
Theorem \ref{result} is obviously smaller than the lattice
$\Lambda_{X}/n$ in Theorem \ref{concini1}. This means that $t_{X}$
makes no sense on the point in $(\Lambda_{X}/n)\backslash
\Lambda_{X}$. The reason comes from the different definitions of
``$\widetilde{X}$" in two toric reduction theorems: Theorem
\ref{result1} and Theorem \ref{concini}. As $t_{X}$ is a function in
$\mathbf{Z}^{s}$, and our discussion in Theorem \ref{result1} is
based on  $\mathbf{R}^{s}$, while  Theorem \ref{concini} is
discussed in $\mathbf{C}^{s}$ (see \cite{Concini}). Hence our
results: Theorems \ref{result1}, \ref{result} are the improvements
of Theorems \ref{concini}, \ref{concini1} respectively.

From Theorem \ref{result}, we see that $t_{X}$ is a sum of some
periodic piecewise polynomials on the shifted sublattice cone of
$\Lambda_{X}^{+}$. This property is shown in the following two
examples.

\begin{exmp}
Let $X=\{1,1,2\}$ be a set of real numbers. Obviously that
$$\mathcal {L} (t_{X})=\frac{1}{(1-e^{-2x})(1-e^{-x})^{2}}.$$
As $1-e^{-2x}=(1+e^{-x})(1-e^{-x})$, we have
$$\hspace{-1cm}\frac{1}{(1-e^{-2x})(1-e^{-x})^{2}}=\frac{(1+e^{-x})^{2}(1-e^{-x})^{2}}{(1-e^{-2x})^{3}(1-e^{-x})^{2}}
=\frac{(1+e^{-x})^{2}}{(1-e^{-2x})^{3}}=(1+e^{-x})^{2}\frac{(e^{2x}D_{x})^{2}}{2!(-2)^{2}}\cdot
\frac{1}{1-e^{-2x}}$$ So  the inverse Laplace transform of
$\frac{1}{(1-e^{-2x})(1-e^{-x})^{2}}$ is equal to
$$t_{X}=(1+\tau_{-1})^{2}\frac{(x+2)(x+4)}{8}t_{2}(x+4)=(1+\tau_{-1})^{2}\frac{(x+2)(x+4)}{8}t_{2}(x)$$
on the lattice cone $\Lambda_{1}^{+}$. Therefore, $t_{X}$  can be
expressed in the following simple form:
$$t_{X}(x)=\left\{
   \begin{aligned}
   &1~~~~~~~~~~~~~~~~~~x=0  \\
   &\frac{(x+2)(x+4)}{4}~~x\in \Lambda_{2}^{+}(1)\\
   &\frac{(x+2)^{2}}{4}~~~~~~~~~~x\in \Lambda_{2}^{+}(2)\\
   \end{aligned}
   \right.$$
\end{exmp}

\begin{exmp}
Let $X=\{({1 \atop 0}),({0 \atop 1}),({-1 \atop ~2})\}$ be a set of
real numbers. Observe that
$$\mathcal {L} (t_{X})=\frac{1}{(1-e^{-x})(1-e^{-y})(1-e^{-(-x+2y)})},$$
we have
\begin{equation*}
\begin{split}&\frac{1}{(1-e^{-x})(1-e^{-y})(1-e^{-(-x+2y)})}\\
=&\frac{(1+e^{-y})(1-e^{-y})-e^{-2y}e^{-x+2y}(1-e^{-(-x+2y)})}{(1-e^{-x})^{2}(1-e^{-y})(1-e^{-(-x+2y)})}\\
=&\frac{1+e^{-y}}{(1-e^{-x})^{2}(1-e^{-(-x+2y)})}-\frac{e^{-2y}e^{-x+2y}}{(1-e^{-x})^{2}(1-e^{-y})}\\
=&(1+e^{-y})\frac{e^{x}D_{(2,1)}}{-2}\cdot\frac{1}{(1-e^{-x})(1-e^{-(-x+2y)})}-(e^{-x})\frac{e^{x}D_{(1,0)}}{-1}
\cdot\frac{1}{(1-e^{-x})(1-e^{-y})},
\end{split}
\end{equation*}
where we use the identity
\begin{eqnarray*}1-e^{-x}&=&(1-e^{-(2y-(-x+2y))})\\
&=&(1-e^{-2y})+e^{-2y}(1-e^{-x+2y})\\
&=&(1+e^{-y})(1-e^{-y})-e^{-2y}e^{-x+2y}(1-e^{-(-x+2y)}).\end{eqnarray*}
Let $A_1=\{({1 \atop 0}),({-1 \atop ~2})\}$ and $A_2=\{({1 \atop
0}),({0 \atop 1})\}$, then the inverse Laplace transform of
$\frac{1}{(1-e^{-x})(1-e^{-y})(1-e^{-(-x+2y)})}$ is equal to
$$t_{X}=(1+\tau_{(0,-1)})\frac{2x+y+2}{2}t_{A_1}(x+1,y)-\tau_{(-1,0)}(x+1)t_{A_2}(x+1,y)$$
 on the lattice cone $\Lambda_{X}^{+}$.
 Therefore,
 $$t_{X}=\frac{2x+y+2}{2}t_{A_1}(x,y)+\frac{2x+y+1}{2}t_{A_1}(x,y-1)-xt_{A_2}(x,y).$$
\end{exmp}

\section*{Acknowledgements}
\hspace{0.5cm}The work was partly supported by the National Natural
Science Foundation of China (Grant Nos. 60533060,10801024, and
U0935004) and the Innovation Foundation of the Key Laboratory of
High-Temperature Gasdynamics of CAS, China.


\end{document}